\documentclass[twoside,10pt]{article}
\usepackage[]{amsmath,amsthm,amssymb}
\usepackage[latin1]{inputenc}

\begin{document}

\begin{center}{{\Large\bf  Integral and series transformations via Ramanujan's} }\end{center}
\begin{center}{{\Large\bf  identities and  Salem's type equivalences } }\end{center}
\begin{center}{{\Large\bf   to the  Riemann hypothesis} }\end{center}
\vspace{0,5cm}
\begin{center}{Semyon  YAKUBOVICH}\end{center}

\markboth{\rm \centerline{ Semyon  YAKUBOVICH}}{}
\markright{\rm \centerline{Integral and Series Transformations}}

\begin{abstract} {\noindent We consider  integral and
series transformations, which are associated \ with Ramanujan's
identities, involving arithmetic functions $a(n), \omega(n),
\sigma_a(n), \ d(n),\  \mu(n),\  \lambda(n), \varphi(n)$ and a ratio
of products of Riemann's zeta functions of different arguments.
Reciprocal inversion formulas are proved in a Banach space of
functions whose Mellin's transforms are integrable over the vertical
line ${\rm Re}\ s > 1$. Examples of new transformations like
Widder-Lambert and Kontorovich-Lebedev type are exhibited.
Particular cases include familiar Lambert and M$\ddot{o}$bius
transformations. Finally a class of equivalences of the Salem type
to the Riemann hypothesis is established.}

\end{abstract}
\vspace{4mm}

\vspace{2mm}


\vspace{4mm}

\section {Introduction and auxiliary results}

Integral and series transformations, which will be derived in the
sequel are based on remarkable Ramanujan's identities involving
arithmetic and Riemann's zeta-functions [6], [10], namely
$$\frac{\zeta(s)\zeta(s-a)\zeta\left(s- b\right)\zeta\left(s -a-b\right)}
{\zeta(2s-a-b)} =\sum_{n=1}^\infty \frac{
\sigma_{a}(n)\sigma_b(n)}{n^{s}},\eqno(1.1)$$
where ${\rm Re}\  s > \max \{1, {\rm Re} \ a +1, {\rm Re}\  b +1,
{\rm Re}\  (a+b) +1\}$,
$$\zeta\left(s\right)\zeta\left(s -a\right)
=\sum_{n=1}^\infty \frac{ \sigma_{a}(n)}{n^{s}},\eqno(1.2)$$
where ${\rm Re} s > \max \{1, {\rm Re}\  a +1\},$
$$\frac{\zeta^2(s)}
{\zeta(2s)} =\sum_{n=1}^\infty \frac{ 2^{\omega(n)}}{n^{s}}, \quad
{\rm Re} s > 1,\eqno(1.3)$$
$$\zeta^2(s) =\sum_{n=1}^\infty \frac{ d(n)}{n^{s}}, \quad
{\rm Re} s > 1,\eqno(1.4)$$
$$\frac{1}{\zeta(s)} =\sum_{n=1}^\infty \frac{\mu(n)}{n^{s}}, \quad {\rm Re} s
> 1,\eqno(1.5)$$
$$\frac{\zeta(s)}
{\zeta(2s)} =\sum_{n=1}^\infty \frac{\left|\mu(n)\right|}{n^{s}},
\quad {\rm Re} s > 1,\eqno(1.6)$$
$$\frac{\zeta(2s)}
{\zeta(s)} =\sum_{n=1}^\infty \frac{\lambda(n)}{n^{s}}, \quad {\rm
Re} s > 1,\eqno(1.7)$$
$$\frac{\zeta^3(s)}
{\zeta(2s)} =\sum_{n=1}^\infty \frac{d\left(n^2\right)}{n^{s}},
\quad {\rm Re} s > 1,\eqno(1.8)$$
$$\frac{\zeta^4(s)}
{\zeta(2s)} =\sum_{n=1}^\infty \frac{d^2(n)}{n^{s}}, \quad {\rm Re}
s > 1,\eqno(1.9)$$
$$\frac{\zeta(s-1)}
{\zeta(s)} =\sum_{n=1}^\infty \frac{\varphi(n)}{n^{s}}, \quad {\rm
Re} s > 2,\eqno(1.10)$$
$$\frac{1- 2^{1-s}}{1-2^{-s}}\zeta(s-1) =\sum_{n=1}^\infty \frac{a(n)}{n^{s}},
\quad {\rm Re} s > 2.\eqno(1.11)$$
Here $\zeta(s)$ is the Riemann zeta-function [10], which satisfies
the familiar functional equation
$$\zeta(s)= 2^s \pi^{s-1} \sin\left({\pi s\over
2}\right)\Gamma(1-s)\zeta(1-s),\eqno(1.12)$$
where $\Gamma(z)$ is Euler's gamma-function, and in the half-plane
${\rm Re}  s= c_0 > 1$ it is represented by the absolutely and
uniformly convergent series with respect to $t \in \mathbb{R},\ s=
c_0+it$
$$\zeta(s) =\sum_{n=1}^\infty \frac{1}{n^{s}},\eqno(1.13)$$
and by the uniformly convergent series
$$(1- 2^{1-s})\zeta(s) =\sum_{n=1}^\infty \frac{(-1)^{n-1}}{n^{s}},\   {\rm Re} s > 0.\eqno(1.14)$$
Further, $a(n)$ in (1.11) denotes the greatest odd divisor of $n$,
$\sigma_{a}(n)$ in (1.1), (1.2) is the sum of $a$-th powers of the
divisors of $n \in \mathbb{N}$. In particular, for pure imaginary
$a=i\tau$ $|\sigma_{i\tau}(n)| \le d(n)$, where $d(n)$ is the
Dirichlet divisor function, i.e. the number of divisors of $n$,
including $1$ and $n$ itself. It has the estimate [10] $d(n)=
O(n^\varepsilon), \ n \to \infty, \ \varepsilon
> 0$. The M$\ddot{o}$bius function is denoted by $\mu(n)$ and
$|\mu(n)| \le 1$. The symbol $\omega(n)$ in (1.3) represents the
number of distinct prime factors of $n$ and it behaves as
$\omega(n)= O(\log\log n), \ n \to \infty$ (see in [8]). By
$\varphi(n)$ Euler's totient function is denoted and  its asymptotic
behavior satisfies [cf. [8]) $\varphi(n)= O\left(n [\log\log
n]^{-1}\right), \ n \to \infty$. Finally, $\lambda(n)$ in (1.7) is
the Liouville function, $|\lambda(n)| \le 1$.

Following similar ideas presented in [12], [14] we define a special
functional space ${\cal M}^{-1}(L_c)$, which will be suitable for
our investigation of the series transformations with arithmetic
functions.

{\bf Definition 1}. Denote by ${\cal M}^{-1}(L_c)$ the space of
functions $f(x), x \in \mathbb{R}_+$,  representable by inverse
Mellin transform of integrable functions $f^{*}(s) \in L_{1}(c)$ on
the vertical line $c =\{s \in \mathbb{C}: {\rm Re} s=c_0 \}$:
$$ f(x) = {1\over 2\pi i} \int_c f^{*}(s)x^{-s}ds.\eqno(1.15)$$

The  space ${\cal M}^{-1}(L_c)$  with  the  usual operations  of
addition   and multiplication by scalar is a linear vector space. If
the norm in ${\cal M}^{-1}(L_c)$ is introduced by the formula
$$ \big\vert\big\vert f \big\vert\big\vert_{{\cal
M}^{-1}(L_c)}= {1\over 2\pi }\int^{+\infty}_{-\infty} | f^{*}(c_0
+it)| dt,\eqno(1.16)$$
then it becomes  a Banach space. Simple properties of the space
${\cal M}^{-1}(L_c)$ follow immediately from Definition 1 and the
basic properties of the Fourier and Mellin transforms of integrable
functions. For instance, the Riemann-Lebesgue lemma yields that
$x^{c_0}f(x)$ is  uniformly bounded, continuous on $\mathbb{R}_+$
and  $x^{c_0}f(x)= o(1)$, when $x\to +\infty $ and $x\to 0$.
Moreover, if $f(x),\  g(x) \in {\cal M}^{-1}(L_c)$, where $g(x)$ is
the inverse Mellin transform (1.14) of the function $g^*(s)$, then
$x^{c_0}f(x)g(x) \in {\cal M}^{-1}(L_c)$ because the product
$x^{c_0}f(x)g(x)$ is the  inverse Mellin transform of the function
$${1\over 2\pi i}\int_c  f^{*}(\tau)g^{*}(s-\tau+ c_0)d\tau,$$
which belongs to $L_1(c)$ by Fubini's theorem. Finally we note that
if $f(x) \in {\cal M}^{-1}(L_c)$  and $x^{c_0-1}g(x) \in
L_1(\mathbb{R}_+)$, then the Mellin convolution
$$\int_0^\infty g(u)f\left({x\over u}\right){du\over u} \in {\cal M}^{-1}(L_c).$$
In fact, the latter integral is  an inverse Mellin transform of the
function $f^{*}(s)g^{*}(s)$ and since $f^{*}(s) \in L_1(c)$ and
$g^{*}(s)$ is essentially bounded on $c$, we have $f^{*}(s)g^{*}(s)
\in L_1(c)$.

A more general space ${\cal M}_{c_1,c_2}^{-1}(L_c)$, which will be
involved as well is defined similarly to the one in [12], [14].

{\bf Definition 2}.  Let $c_1, c_2 \in \mathbb{R}$ be such that $2
\hbox{sign}\ c_1 + \hbox{sign}\  c_2 \ge 0$. By ${\cal
M}_{c_1,c_2}^{-1}(L_c)$ we denote the space of functions $f(x), x
\in \mathbb{R}_+$, representable in the form (1.15), where
$s^{c_2}e^{\pi c_1|s|} f^*(s) \in L_1(c)$.

It is a Banach space with the norm
$$ \big\vert\big\vert f \big\vert\big\vert_{{\cal
M}_{c_1,c_2}^{-1}(L_c)}= {1\over 2\pi }\int_{c} e^{\pi c_1|s|}
|s^{c_2} f^{*}(s) ds|, \ {\rm Re} s = c_0.$$,

\section{Transformations with arithmetic functions}

We begin with the following result.

{\bf Theorem 1}. {\it Let $f \in {\cal M}^{-1}(L_c), \ c_0 > 1$.
Then for all $x >0$ the following series expansions with the
M$\ddot{o}$bius function are true
$$f(x)= \sum_{n=1}^\infty \mu(n)\sum_{m=1}^\infty
f(xnm),\eqno(2.1)$$
$$f(x)= \sum_{n=1}^\infty \sum_{m=1}^\infty \mu(m)
f(xnm),\eqno(2.2)$$
$$f(x)=\sum_{k=0}^\infty \sum_{n=1}^\infty 2^k \mu(n)\sum_{m=1}^\infty (-1)^{m-1}
f\left(xnm 2^k\right),\eqno(2.3)$$
$$f(x)=\sum_{n=1}^\infty (-1)^{n-1} \sum_{k=0}^\infty \sum_{m=1}^\infty 2^k \mu(m)
f\left(xnm 2^k\right),\eqno(2.4)$$
$$f(x)- 2f(2x)= \sum_{n=1}^\infty \mu(n)\sum_{m=1}^\infty
(-1)^{m-1} f(xnm).\eqno(2.5)$$
Moreover, expansions $(2.1), (2.2)$ generate reciprocal pair of
transformations
$$g(x)= {1\over 2\pi i}\int_c  \zeta(s)f^{*}(s) x^{-s}ds= \sum_{n=1}^\infty
f(xn),\eqno(2.6)$$
$$f(x)= {1\over 2\pi i}\int_c  \frac{g^{*}(s)}{\zeta(s)} x^{-s}ds=
\sum_{n=1}^\infty \mu(n) g(xn),\eqno(2.7)$$
which are automorphisms of the space ${\cal M}^{-1}(L_c)$ and
satisfy the following inequalities for the norms
$$[\zeta(c_0)]^{-1}\ ||g||_{{\cal M}^{-1}(L_c)} \le ||f||_{{\cal
M}^{-1}(L_c)}\le \zeta(c_0)\ ||g||_{{\cal M}^{-1}(L_c)}, \ c_0 >
1.\eqno(2.8)$$
Analogously expansions  $(2.3), (2.4)$ generate reciprocal
transformations
$$g(x)= {1\over 2\pi i}\int_c  (1-2^{1-s}) \zeta(s)f^{*}(s) x^{-s}ds= \sum_{n=1}^\infty
(-1)^{n-1} f(xn),\eqno(2.9)$$
$$f(x)= {1\over 2\pi i}\int_c  \frac{g^{*}(s)}{(1-2^{1-s})\zeta(s)} x^{-s}ds=
\sum_{k=0}^\infty\sum_{n=1}^\infty 2^k \mu(n) g\left(xn
2^k\right),\eqno(2.10)$$
which are automorphisms of the space ${\cal M}^{-1}(L_c)$ and
satisfy the norm estimates}
$$ [\zeta(c_0)]^{-1} ||g||_{{\cal M}^{-1}(L_c)} \le ||f||_{{\cal
M}^{-1}(L_c)}\le (1- 2^{1-c_0})^{-1} \zeta(c_0) ||g||_{{\cal
M}^{-1}(L_c)}, \ c_0 > 1.\eqno(2.11)$$

\begin{proof}   In fact, the validity of equalities (2.1)- (2.5)
follows immediately from  representation (1.15), identities (1.5),
(1.13), (1.14) and elementary sum of geometric progression after the
change of the order of summation and integration via Fubini's
theorem owing to the  following estimates
$$\sum_{n=1}^\infty |\mu(n)|\sum_{m=1}^\infty |f(xnm)| \le
{x^{-c_0}\zeta^2(c_0)\over 2\pi}  \int_c |f^*(s) ds| < \infty, \ x
>0,$$
$$\sum_{k=0}^\infty\sum_{n=1}^\infty 2^k |\mu(n)|\sum_{m=1}^\infty |f\left(xnm
2^k\right)| \le {x^{-c_0}\zeta^2(c_0)\over 2\pi (1- 2^{1-c_0})}
\int_c |f^*(s) ds| < \infty, \ x >0.$$
Hence we establish reciprocal equalities (2.6), (2.7), where
$f^*(s)= g^*(s)[\zeta(s)]^{-1}$ by virtue of the uniqueness theorem
for the Mellin transform. This also guarantees the automorphism of
the space ${\cal M}^{-1}(L_c)$ under transformations (2.6), (2.7).
Finally, since (see (1.16))
$$||f||_{{\cal M}^{-1}(L_c)} = {1\over 2\pi} \int_c \left|\zeta(s) f^*(s) {ds\over \zeta(s)}\right |
\le \zeta(c_0)\ ||g||_{{\cal M}^{-1}(L_c)},$$
$$||g||_{{\cal M}^{-1}(L_c)} =  {1\over 2\pi} \int_c \left|\zeta(s) f^*(s) ds\right |
\le \zeta(c_0)\ ||f||_{{\cal M}^{-1}(L_c)},$$
we prove inequalities (2.8). In the same manner we establish the
automorphism of the space ${\cal M}^{-1}(L_c)$ under  reciprocal
pair (2.9), (2.10) and estimates ($g^*(s)=
f^*(s)\zeta(s)(1-2^{1-s})$)
$$||g||_{{\cal M}^{-1}(L_c)} = {1\over 2\pi} \int_c \left|(1-2^{1-s})\zeta(s) f^*(s) ds\right |
\le \zeta(c_0)\ ||f||_{{\cal M}^{-1}(L_c)},$$
$$||f||_{{\cal M}^{-1}(L_c)} =  {1\over 2\pi} \int_c \left|{g^*(s) ds\over (1-2^{1-s})\zeta(s)}\right |
\le (1-2^{1-c_0})^{-1} \zeta(c_0)\ ||g||_{{\cal M}^{-1}(L_c)},$$
yield inequalities (2.11).
\end{proof}

Analogously, calling Ramanujan's identities (1.2)- (1.10) we come
out with two more theorems, which we leave without proof.

{\bf Theorem 2}. {\it Let $f \in {\cal M}^{-1}(L_c), \ c_0 > 1$.
Then for all $x >0$ the following series expansions with arithmetic
functions hold valid}
$$f(x)= \sum_{n=1}^\infty |\mu(n)|\sum_{m=1}^\infty \lambda (m)
f(xnm),$$
$$f(x)= \sum_{n=1}^\infty \lambda (n) \sum_{m=1}^\infty |\mu(m)|
f(xnm),$$
$$f(x)=\sum_{k,n=1}^\infty \mu(k) \lambda(n)\sum_{m=1}^\infty
2^{\omega(m)} f\left(xnmk\right),$$
$$f(x)=\sum_{n=1}^\infty 2^{\omega(n)} \sum_{k,m=1}^\infty \lambda(k)\mu(m)
f\left(xnmk\right),$$
$$f(x)= \sum_{k,m=1}^\infty \mu(k)\mu(m)\sum_{n=1}^\infty
d(n) f(xnmk),$$
$$f(x)= \sum_{m=1}^\infty d(m)\sum_{k,n=1}^\infty
\mu(k)\mu(n) f(xnmk),$$
$$f(x)= \sum_{k,n,j=1}^\infty \mu(k)\mu(j)\lambda(n)\sum_{m=1}^\infty
d(m^2) f(xnmjk),$$
$$f(x)= \sum_{m=1}^\infty d(m^2)\sum_{k,n,j=1}^\infty
\mu(k)\mu(n) \lambda(j) f(xnmkj),$$
$$f(x)= \sum_{k,n,j,l=1}^\infty \mu(k)\mu(j)\mu(l)\lambda(n)\sum_{m=1}^\infty
d^2(m) f(xnmjkl),$$
$$f(x)= \sum_{m=1}^\infty d^2(m)\sum_{k,n,j,l=1}^\infty
\mu(k)\mu(n)\mu(l) \lambda(j) f(xnmkjl),$$
$$f(x)= \sum_{k,n=1}^\infty k \mu(k)\sum_{m=1}^\infty
\varphi(m) f(xnmk), \ c_0 > 2,$$
$$f(x)= \sum_{m=1}^\infty \varphi(m)\sum_{k,n=1}^\infty k
\mu(k) f(xnmk),\ c_0 > 2,$$
$$f(x)= \sum_{k,n=1}^\infty n^{a} \mu(k)\mu(n)\sum_{m=1}^\infty
\sigma_a(m) f(xnmk), \ c_0 > \max \{1, {\rm Re}\  a +1\},$$
$$f(x)= \sum_{m=1}^\infty \sigma_a(m)\sum_{k,n=1}^\infty k^a
\mu(k) \mu(n) f(xnmk),\ c_0 > \max \{1, {\rm Re}\  a +1\}.$$

{\bf Theorem 3}. {\it Let $f \in {\cal M}^{-1}(L_c),\ c_0 > 1$. Then
for all $x >0$ the following reciprocal series transformations are
automorphisms in ${\cal M}^{-1}(L_c)$ with the corresponding norm
estimates, namely}
$$g(x)=\sum_{m=1}^\infty \lambda (m) f(xm),$$
$$f(x)= \sum_{n=1}^\infty |\mu(n)| g(xn),$$
$$[\zeta(c_0)]^{-1} \ ||g||_{{\cal M}^{-1}(L_c)} \le ||f||_{{\cal
M}^{-1}(L_c)}\le {\zeta( c_0)\over \zeta(2 c_0)}\ ||g||_{{\cal
M}^{-1}(L_c)}, \ c_0 > 1;$$
$$g(x)=\sum_{n=1}^\infty 2^{\omega(n)} f(xn),$$
$$f(x)= \sum_{k,m=1}^\infty \lambda(k)\mu(m) g\left(xmk\right),$$
$${\zeta(2c_0)\over \zeta^2(c_0)}\ ||g||_{{\cal M}^{-1}(L_c)} \le ||f||_{{\cal
M}^{-1}(L_c)}\le \zeta^2(c_0)\ ||g||_{{\cal M}^{-1}(L_c)}, \ c_0
> 1;$$
$$g(x)= \sum_{n=1}^\infty d(n) f(xn),$$
$$f(x)= \sum_{k,n=1}^\infty \mu(k)\mu(n)g(xkn),$$
$$\zeta^{-2}(c_0)\ ||g||_{{\cal M}^{-1}(L_c)} \le ||f||_{{\cal
M}^{-1}(L_c)}\le \zeta^2(c_0)\ ||g||_{{\cal M}^{-1}(L_c)}, \ c_0
> 1;$$
$$g(x)= \sum_{m=1}^\infty d(m^2) f(xm),$$
$$f(x)= \sum_{k,n,j=1}^\infty \mu(k)\mu(j)\lambda(n)g(xnjk),$$
$$\zeta^{-3}(c_0)\ ||g||_{{\cal M}^{-1}(L_c)} \le ||f||_{{\cal
M}^{-1}(L_c)}\le \zeta^3(c_0)\zeta(2c_0)\ ||g||_{{\cal
M}^{-1}(L_c)}, \ c_0 > 1;$$
$$g(x)= \sum_{m=1}^\infty d^2(m) f(xm),$$
$$f(x)= \sum_{k,n,j,l=1}^\infty \mu(k)\mu(j)\mu(l)\lambda(n)g(xnjkl),$$
$$\zeta^{-4}(c_0)\ ||g||_{{\cal M}^{-1}(L_c)} \le ||f||_{{\cal
M}^{-1}(L_c)}\le \zeta^4(c_0)\zeta(2c_0)\ ||g||_{{\cal
M}^{-1}(L_c)}, \ c_0 > 1;$$
$$g(x)= \sum_{m=1}^\infty \varphi(m) f(xm),$$
$$f(x)= \sum_{k,n=1}^\infty k \mu(k)g(xnk),$$
$$[\zeta(c_0)\zeta(c_0-1)]^{-1}\ ||g||_{{\cal M}^{-1}(L_c)} \le ||f||_{{\cal
M}^{-1}(L_c)}\le \zeta(c_0)\zeta(c_0-1)\ ||g||_{{\cal M}^{-1}(L_c)},
\ c_0 > 2;$$
$$g(x)= \sum_{m=1}^\infty \sigma_a(m) f(xm),$$
$$f(x)= \sum_{k,n=1}^\infty n^{a} \mu(k)\mu(n)g(xnk),$$
$$[\zeta(c_0)\zeta(c_0-{\rm Re} a)]^{-1}\ ||g||_{{\cal M}^{-1}(L_c)} \le ||f||_{{\cal
M}^{-1}(L_c)}\le \zeta(c_0)\zeta(c_0-{\rm Re} a)\ ||g||_{{\cal
M}^{-1}(L_c)}, \ c_0 > \max \{1, {\rm Re}\  a +1\}.$$

 Let $f(x)= e^{-x}$, which  evidently belongs to the
space ${\cal M}^{-1}(L_c)$ since $f^*(s)=\Gamma(s)$ is the Euler
gamma-function. Substituting it in (2.1)-(2.5) and calculating
elementary series we come out to the Lambert type expansions  (cf.
[8])
$$e^{-x}= \sum_{n=1}^\infty {\mu(n)\over e^{xn}-1},\  x >0,$$
$$e^{-x}=\sum_{k=0}^\infty \sum_{n=1}^\infty {2^k \mu(n)\over \exp\left(xn 2^k\right) +1}, \ x >0,$$
$$e^{-x}- 2e^{-2x}= \sum_{n=1}^\infty {\mu(n)\over e^{xn} +1}, \ x >0.$$
Further, the Parseval equality for the Mellin transform [11] and
Fubini's theorem allow to write the modified Laplace transform [3]
of $ f \in {\cal M}^{-1}(L_c)$ in the form
$$\int_0^\infty e^{-x/t} f(t){dt\over t}= {1\over 2\pi i}\int_c
\Gamma(s)f^*(s) x^{-s} ds.\eqno(2.12)$$
Moreover, due to Definition 2 and Stirling's asymptotic formula for
gamma-functions [1, Vol. 1] it forms a bijective map of the space
${\cal M}^{-1}(L_c)$ onto its subspace ${\cal M}_{1/2,1/2
-c_0}^{-1}(L_c)$. Thus appealing to Theorem 1 we will derive the
Widder type inversion formulas for the Lambert transform (see in
[2], [4], [13]) and Widder-Lambert  type transforms. Precisely, we
prove

{\bf Theorem 4}. {\it Let $f \in {\cal M}^{-1}(L_c)$ and $c_0 >1$.
Then the modified  Lambert transform
$$g(x)= \int_0^\infty {f(t) dt\over t(e^{x/t} -1)}, \ x > 0\eqno(2.13)$$
maps bijectively onto the space ${\cal M}_{1/2,1/2 -c_0}^{-1}(L_c)$
and for all $x >0$ the Widder type inversion formula holds true
$$f(x)= \lim_{k\to \infty} \left(-x{d\over dx} \right)\prod_{j=1}^k
\left(1- {x\over j}{d\over dx}\right)\sum_{n=1}^\infty \mu(n)
g(xkn).\eqno(2.14)$$
Analogously, the Widder-Lambert type transformation

$$g(x)= \int_0^\infty {f(t) dt\over t(e^{x/t} +1)}, \ x > 0\eqno(2.15)$$
is a bijective map between spaces ${\cal M}^{-1}(L_c)$,  ${\cal
M}_{1/2,1/2 -c_0}^{-1}(L_c)$ and for all $x >0$ the following
inversion formula takes place}
$$f(x)= \lim_{k\to \infty} \left(-x{d\over dx} \right)\prod_{j=1}^k
\left(1- {x\over j}{d\over dx}\right)
\sum_{m=0}^\infty\sum_{n=1}^\infty 2^m \mu(n) g\left(xkn
2^m\right).\eqno(2.16)$$

\begin{proof} In fact, the proof is based on Theorem 1, equality
(2.12), a familiar infinite product for the gamma-function (see, for
instance, [14, p.48])
$${1\over \Gamma(s)}= \lim_{k\to \infty} s k^{-s}\prod_{j=1}^k
\left(1+ {s\over j}\right),$$
and the asymptotic behavior $|\Gamma(s)|^{-1}=
e^{\pi|s|/2}|s|^{1/2-c_0}, \ s= c_0+it, \ |t| \to \infty$ via
Stirling formula. So owing to Theorem 1 and the absolute and uniform
convergence, which guarantees the change of the order of integration
and summation,  the modified Lambert transform bijectively maps
${\cal M}^{-1}(L_c)$ onto  ${\cal M}_{1/2,1/2 -c_0}^{-1}(L_c)$ and
represented by (2.13), namely
$$g(x)={1\over 2\pi i}\int_c \zeta(s)\Gamma(s)f^*(s) x^{-s} ds
= \sum_{n=1}^\infty \int_0^\infty e^{-xn/t} f(t){dt\over t}$$
$$= \int_0^\infty {f(t) dt\over t(e^{x/t} -1)}, x >0.$$
Reciprocally, following similarly to [14, p.49] and appealing to the
Lebesgue dominated convergence theorem and equality (2.7),  we find
$$f(x)= {1\over 2\pi i}\int_c {g^*(s) x^{-s}\over \zeta(s)\Gamma(s)} ds
= \lim_{k\to \infty} {1\over 2\pi i}\int_c \prod_{j=1}^k \left(1+
{s\over j}\right) {s g^*(s) (kx)^{-s}\over \zeta(s)}
ds$$
$$= \lim_{k\to \infty} \left(-x{d\over dx} \right)\prod_{j=1}^k \left(1- {x\over j}{d\over
dx}\right)\sum_{n=1}^\infty \mu(n) g(xkn),$$
which gives (2.14). In the same manner, employing again (1.14) and
Theorem 1 we deduce the representation (2.15) of the Widder-Lambert
type transform
$$g(x)={1\over 2\pi i}\int_c (1-2^{1-s})\zeta(s)\Gamma(s)f^*(s) x^{-s} ds
= \sum_{n=1}^\infty (-1)^{n-1} \int_0^\infty e^{-xn/t} f(t){dt\over
t}$$
$$= \int_0^\infty {f(t) dt\over t(e^{x/t} +1)}, x >0.$$
Finally, the same motivations perform the chain of equalities
$$f(x)= {1\over 2\pi i}\int_c {g^*(s) x^{-s}\over (1-2^{1-s})\zeta(s)\Gamma(s)} ds
= \lim_{k\to \infty} {1\over 2\pi i}\sum_{m=0} 2^m \int_c
\prod_{j=1}^k \left(1+ {s\over j}\right) {s g^*(s) \left(kx
2^m\right)^{-s}\over \zeta(s)} ds$$
$$= \lim_{k\to \infty} \left(-x{d\over dx} \right)\prod_{j=1}^k \left(1- {x\over j}{d\over
dx}\right)\sum_{m=0}^\infty \sum_{n=1}^\infty 2^m \mu(n) g\left(xkn
2^m\right),$$
which, in turn,  yield (2.16).
\end{proof}

Transformation (2.15) can be generalized considering the following
two-parametric family of functions
$$U_{k,m}(x)={1\over 2\pi i}\int_c
[(1-2^{1-s})\zeta(s)]^{k+1}\Gamma^{m+1}(s)x^{-s} ds, x >0, \ k, m
\in \mathbb{N}_0.\eqno(2.17)$$
The case $k=m$ we denote by $U_k(x)$. The case $k=m=0$ gives
$U_0(x)= (e^x+1)^{-1}$. One can express the kernel (2.17) in terms
of the iterated Mellin convolution. Indeed, via (1.14) and simple
calculations we obtain
$$U_{k,m}(x)= \sum_{n_1,n_2,\dots,n_{k-m}=1}^\infty
(-1)^{\sum_{j=1}^{k-m} n_j- k+m}$$ $$ \times \int_{\mathbb{R}^{m}_+}
\prod_{j=1}^m (e^{u_j}+1)^{-1} \left(\exp\left({x n_1 n_2\dots
n_{k-m}\over u_1 u_2 \dots u_m}\right) + 1\right)^{-1} {du_1 du_2
\dots du_m \over u_1 u_2 \dots u_m},\   k > m,\eqno(2.18)$$
$$U_{k,m}(x)\equiv U_k(x) = \int_{\mathbb{R}^{k}_+}
\left(\exp\left({x \over u_1 u_2 \dots u_k}\right) + 1\right)^{-1}
\prod_{j=1}^{k} (e^{u_j}+1)^{-1}\  {du_j \over u_j},\quad
k=m,\eqno(2.19)$$
$$ U_{k,m}(x) = \int_{\mathbb{R}^{m}_+} \prod_{j=1}^{k+1}
(e^{u_j}+1)^{-1}\exp\left(- \sum_{j=k+2}^{m} u_j\right) \
\exp\left(- {x \over u_1 u_2 \dots u_m}\right) {du_1\dots du_m \over
u_1 u_2 \dots u_m},\quad k< m.\eqno(2.20)$$

Thus an analog of Theorem 4 will be

{\bf Theorem 5} {\it  Let $f \in {\cal M}^{-1}(L_c)$ and $c_0
>1$. Then the integral transformation
$$g(x)= \int_0^\infty U_{k,m}\left({x\over t}\right)f(t) {dt\over t},
\ x > 0\eqno(2.21)$$
is a bijective map between spaces ${\cal M}^{-1}(L_c)$,  ${\cal
M}_{(m+1)/2,(m+1)(1/2 -c_0)}^{-1}(L_c)$ and for all $x >0$ the
following inversion formula takes place}
\begin{equation*}
\begin{split}
f(x)= \lim_{l\to \infty} \left(-x{d\over dx} \right)^{m+1}
\prod_{j=1}^l \left(1- {x\over j}{d\over dx}\right)^{m+1}\\
\times\sum_{j_1,\dots,j_k =0}^\infty\sum_{n_1,\dots,n_k =1}^\infty
\prod_{i=1}^k 2^{j_i}\mu(n_i) g\left(xl^{m+1} \prod_{i=1}^k 2^{j_i}
n_i \right).
\end{split}
\end{equation*}

\section{Transformations of the Kontorovich-Lebedev type}

The familiar Kontorovich-Lebedev transform (see for instance in [9],
[14], [16]) is defined by
$$K_{i\tau}[f]= \int_0^\infty K_{i\tau}(x) f(x)dx,\ \tau \in \mathbb{R}_+,\eqno(3.1)$$
where the integral converges in an appropriate sense and
$K_{\nu}(x), \ \nu \in \mathbb{C}, x >0$ is the modified Bessel
function [1, Vol. 2] having the following integral representations
$$K_{\nu}(2\sqrt x)= {1\over 4\pi i} \int_{a- i\infty}^{a+
i\infty} \Gamma\left(s+ {\nu\over 2}\right)\Gamma\left(s- {\nu\over
2}\right)x^{-s}ds, \ a > |{\rm Re } \nu|,\eqno(3.2)$$
$$K_\nu(x)= \int_0^\infty e^{-  x\cosh u}\cosh \nu u \ du.\eqno(3.3)$$
The main goal of this section is to consider an analog of the
Kontorovich-Lebedev transform (3.1) involving the kernel, which we
will call the  Macdonald-Lambert function ${\cal M}_{\nu}(x)$,
represented by
$${\cal M}_{\nu}(x)= \int_0^\infty {\cosh \nu u \ du\over e^{ x\cosh u} -
1}, \ x >0, \ \nu \in \mathbb{C}.\eqno(3.4)$$
Precisely, letting in (3.4) $\nu$ as a pure imaginary number,
$\nu=i\tau, \ \tau >0$ let us consider the following transformation
$${\cal M}_{i\tau}[f]= \int_0^\infty {\cal M}_{i\tau}(x) f(x)dx,\ \tau \in \mathbb{R}_+.\eqno(3.5)$$
First we observe via (3.4),  that ${\cal M}_{i\tau}(x)$ is a
real-valued function. Moreover, due to (3.3) and elementary
summation it can be represented by the following series of the
modified Bessel functions
$${\cal M}_{i\tau}(x)= \sum_{n=1}^\infty  K_{i\tau}(n x), x >0,\eqno(3.6)$$
where the corresponding change of the order of integration and
summation is by virtue of the absolute and uniform convergence.
Hence invoking the uniform inequality for the modified Bessel
function [16]
$$|K_{i\tau}(x)| \le e^{-r\tau}K_0(x\cos r), \ r \in [0, \pi/2),$$
we have, accordingly, the estimate
$$|{\cal M}_{i\tau}(x)| \le  e^{-r\tau}\sum_{n=1}^\infty  K_{0}(n x\cos r)= e^{-r\tau}{\cal M}_{0}(x\cos r).\eqno(3.7)$$
Meanwhile, making a simple change of variable and shifting the
vertical contour in (3.2) to the right into  the half-plane ${\rm
Re} s > 1$, we substitute it in (3.6). Then inverting  again the
order of integration and summation owing to the absolute and uniform
convergence, we employ (1.13) to deduce the formula
$${\cal M}_{i\tau}(x)= {1\over 2\pi i} \int_{a- i\infty}^{a+
i\infty} 2^{s-2}\Gamma\left({s+ i\tau\over 2}\right)\Gamma\left({s-
i\tau\over 2}\right)\zeta(s) x^{-s}ds,\ a > 1.\eqno(3.8)$$
Reciprocally, taking the Mellin transform of the kernel ${\cal
M}_{i\tau}(x)$,  it yields
$$\int_0^\infty {\cal M}_{i\tau}(x) x^{s-1} dx = 2^{s-2}\Gamma\left({s+ i\tau\over 2}\right)\Gamma\left({s-
i\tau\over 2}\right)\zeta(s),\ {\rm Re}  s= c_0 >1\eqno(3.9)$$
and one can justify the absolute convergence of the integral in
(3.9) shifting a contour in (3.8) to the left and  to the right from
the line ${\rm Re} s = c_0$ in order to get the corresponding
behavior near zero and  infinity, respectively.

 A relationship of (3.5) with  the Kontorovich-Lebedev
transform (3.1) is given by

{\bf Lemma 1}. {\it Let $f \in {\cal M}^{-1}(L_c), c_0=1-a,\ a >1$.
Then for all $\tau \in \mathbb{R}_+$
$${\cal M}_{i\tau}[f]= K_{i\tau}[g],\eqno(3.10)$$ where $g(x)$ is the series
transformation (see $(2.6)$)
$$g(x)={1\over 2\pi i}\int_{a- i\infty}^{a+
i\infty}\zeta(s)f^*(1-s)x^{s-1} ds= \sum_{n=1}^\infty {1\over
n}f\left ({x\over n}\right)\eqno(3.11)$$
and the following equality holds}
$${\cal M}_{i\tau}[f]= {1\over 2\pi i} \int_{a- i\infty}^{a+
i\infty} 2^{s-2}\Gamma\left({s+ i\tau\over 2}\right)\Gamma\left({s-
i\tau\over 2}\right)\zeta(s)f^*(1-s)ds.\eqno(3.12)$$
\begin{proof} In fact, since via conditions of the theorem
$$\int_0^\infty \left|{\cal M}_{i\tau}(x)\right| x^{a-1} \int_{a- i\infty}^{a+
i\infty} |f^*(1-s)ds| dx < \infty,$$
the proof of (3.12) is straightforward by substitution (3.9) into
the right-hand side of (3.12) and inversion of the order of
integration with the use of Fubini's theorem and (3.5). In the same
manner we prove composition (3.10), where $g(x)$ can be represented
by (3.11) similarly to (2.6).
\end{proof}

The main result of this section is an inversion theorem for the
Kontorovich-Lebedev like transformation (3.5). For a different class
of such index transformations and their inversions we refer to [15].
Our method will be based on Sneddon's operational approach to invert
the Kontorovich-Lebedev transform (3.1) (see [9], Ch. 6).

We have

{\bf Theorem 6}. {\it Let $f \in {\cal M}^{-1}(L_c), c_0=1-a,\ a
>1$. Let $f^*(s)$ be  analytic in the strip ${\rm Re} s \in [1-a, 1+a],\ a >1,\ f^*(0)=0$
and $\zeta(-c_0 -it)f^*(1+c_0+it) \in L_1(\mathbb{R}) \cap
L_p(\mathbb{R}), p >1$ for all $c_0 \in [-a, a]$.

If ${\cal M}_{i\tau}[f] \in L_1(\mathbb{R}_+; \tau
e^{\pi\tau}d\tau)$, then for almost all $x \in \mathbb{R}_+$ the
following inversion formula holds
$$x f(x)=  \int_0^\infty \tau \sinh
\pi\tau \ \hat{{\cal M}}_{i\tau}(x){\cal M}_{i\tau}[f] d\tau,\ x
>0,\eqno(3.13)$$
where}
$$\hat{{\cal M}}_{i\tau}(x)= {2^{3/2}\over \sqrt x \
\tau\sinh(\pi\tau/2)}\sum_{n=1}^\infty {\mu(n)\over n^{3/2}}\left(1+
{4\pi^2\over x^2n^2}\right)^{1/4} P_{(i\tau-1)/2}^{1/2}\left(1+
{8\pi^2\over x^2n^2}\right).\eqno(3.14)$$
\begin{proof} We begin substituting integral representation (3.4) with $\nu=i\tau$
into (3.5) and changing the order of integration by Fubini's
theorem. It is indeed allowed via conditions of Lemma 1 and the
estimate
$$ \int_0^\infty  du \int_0^\infty {|f(x)|dx\over e^{ x\cosh u} -
1} \le  {1\over 2\pi}\int_0^\infty  du \int_0^\infty {x^{a-1}\over
e^{ x\cosh u} - 1}dx \int_{a- i\infty}^{a+ i\infty}|f^*(1-s)ds|$$
$$= {1\over 2\pi}\int_0^\infty  {du\over \cosh^a u} \int_0^\infty {x^{a-1}\over
e^{ x} - 1}dx \int_{a- i\infty}^{a+ i\infty}|f^*(1-s)ds| < \infty, \
a > 1.$$
Consequently,
$${\cal M}_{i\tau}[f]= \int_0^\infty \cos\tau u \int_0^\infty {f(x)\over e^{ x\cosh u} -
1}dx du.\eqno(3.15)$$
Hence as we see in the above estimate the inner integral with
respect to $x$ is an integrable function by $u$. Moreover,
inequality (3.7) and conditions of the theorem guarantee that ${\cal
M}_{i\tau}[f] \in L_1(\mathbb{R}_+)$. Thus inverting the cosine
Fourier transform in (3.15) we arrive at the equality
$${2\over \pi} \int_0^\infty {\cal M}_{i\tau}[f] \cos\tau u \ d\tau = \int_0^\infty {f(x) dx \over e^{ x\cosh u} -
1},$$
or after simple substitution $v= \cosh u$ it becomes
$${2\over \pi} \int_0^\infty {\cal M}_{i\tau}[f] \cos\left(\tau \log\left(v+ \sqrt {v^2-1}\right)\right) d\tau = \int_0^\infty {f(x)\over e^{ xv} -
1}dx,\ v >1.\eqno(3.16)$$
A differentiation under integral sign in (3.16) with respect to $v$
is still allowed by virtue of the absolute and uniform convergence
of the corresponding integrals. Precisely, in its left hand-side it
is owing to inequality (3.7) and in the right- hand side by the
inequality
$$\int_0^\infty {x |f(x)| e^{vx}\over (e^{ xv} - 1)^2} dx \le
{1\over 2\pi}\int_0^\infty {x^{a} e^x \over (e^{ x} - 1)^2}dx
\int_{a- i\infty}^{a+ i\infty}|f^*(1-s)ds| < \infty, \ a > 1.$$
Thus (3.16) yields
$${2\over \pi} \int_0^\infty \tau {\cal M}_{i\tau}[f] {\sin\left(\tau \log\left(v+ \sqrt {v^2-1}\right)\right)
\over \sqrt {v^2-1}} d\tau = \int_0^\infty { xf(x) e^{xv} \over (e^{
xv} - 1)^2}dx$$
$$= \int_0^\infty  xf(x) \sum_{n=1}^\infty n e^{- xvn} dx.$$
But
$$\int_0^\infty  x|f(x)| \sum_{n=1}^\infty n e^{- xvn} dx \le
{\zeta(a)\Gamma(a+1)\over 2\pi v^{a+1}} \int_{a- i\infty}^{a+
i\infty}|f^*(1-s)ds|, \ a > 1.$$
Therefore, an interchange of the order of integration and summation
in the right-hand side of latter equality is allowed and after a
simple change of variables we come out with
$${2\over \pi} \int_0^\infty \tau {\cal M}_{i\tau}[f] {\sin\left(\tau \log\left(v+ \sqrt {v^2-1}\right)\right)
\over \sqrt {v^2-1}} d\tau = \int_0^\infty  e^{-vx} x
\sum_{n=1}^\infty {1\over n} f\left({x\over n}\right)
dx.\eqno(3.17)$$
Meanwhile, appealing to the identity (cf. [9, p. 359])
$$ {\sin\left(\tau \log\left(v+ \sqrt {v^2-1}\right)\right)
\over \sqrt {v^2-1}}= {1\over \pi} \sinh \pi\tau \int_0^\infty
e^{-vx}K_{i\tau}(x) dx,$$
and via condition  of the theorem ${\cal M}_{i\tau}[f] \in
L_1(\mathbb{R}_+; \tau e^{\pi\tau}d\tau)$, we substitute the latter
integral into (3.17) and change the order of integration. Then
canceling the Laplace transform due to the uniqueness theorem for
Laplace transform of integrable functions [11],  we arrive for
almost all $x >0$ at the equality
$${2\over \pi^2} \int_0^\infty \tau  \sinh \pi\tau \ K_{i\tau}(x) {\cal M}_{i\tau}[f] d\tau =
\sum_{n=1}^\infty {x\over n} f\left({x\over n}\right).\eqno(3.18)$$
However, the right-hand side of (3.18) is given by the integral
(3.11), which becomes after a simple change of variables as
$$\sum_{n=1}^\infty {x\over n} f\left({x\over n}\right)= {1\over 2\pi i}\int_{- a- i\infty}^{- a+
i\infty}\zeta(-s)f^*(1+s)x^{-s} ds.$$
Moreover,  conditions of the theorem allow to shift the contour to
the right by the Cauchy theorem and write for each $x >0$
$${1\over 2\pi i}\int_{- a- i\infty}^{- a+
i\infty}\zeta(-s)f^*(1+s)x^{-s} ds= {1\over 2\pi i}\int_{a-
i\infty}^{a+ i\infty}\zeta(-s)f^*(1+s)x^{-s} ds.\eqno(3.19)$$
Indeed, since $f^*(0)=0$ and $f^*(s)$ is analytic in the strip ${\rm
Re} s \in [1-a, 1+a],\ a >1,$ we have that the limit of the product
$\zeta(-s)f^*(1+s)$ exists when $s \to -1$ and  it is analytic in
the strip ${\rm Re} s \in [-a, a]$.   Further as in [11], p. 125 the
condition $\zeta(-c_0 -it)f^*(1+c_0+it) \in L_p(\mathbb{R}), p
>1$ for all $c_0 \in [-a, a]$ implies that $\zeta(-s)f^*(1+s)x^{-s},\ s=
c_0+it$ goes to zero when $|t| \to \infty$ uniformly for
$-a+\varepsilon \le c_0 \le a-\varepsilon$ for any small fixed
positive $\varepsilon$. Therefore (3.19) holds. Returning to (3.18)
and accounting (3.2) with Fubini's theorem, which is applicable
under integrability condition on ${\cal M}_{i\tau}[f]$, it becomes
$$ {1\over 2\pi i} \int_{a- i\infty}^{a+
i\infty} 2^{s-2}x^{-s} \left[{2\over \pi^2} \int_0^\infty \tau \sinh
\pi\tau \ \Gamma\left({s+ i\tau\over 2}\right)\Gamma\left({s-
i\tau\over 2}\right) {\cal M}_{i\tau}[f] d\tau\right]ds$$$$ =
{1\over 2\pi i}\int_{a- i\infty}^{a+ i\infty}\zeta(-s)f^*(1+s)x^{-s}
ds.\eqno(3.20)$$
Canceling the inverse Mellin transform from both sides of (3.20),
because the integrands are $L_1$-functions and dividing by
$\zeta(-s)$, we obtain
$$ f^*(1+s)= {2^{s-1}\over \zeta(-s)\pi^2} \int_0^\infty \tau \sinh
\pi\tau \ \Gamma\left({s+ i\tau\over 2}\right)\Gamma\left({s-
i\tau\over 2}\right) {\cal M}_{i\tau}[f] d\tau.$$
Hence taking the inverse Mellin transform over $(b-i\infty,
b+i\infty),\ 1< b< 2$ from both sides of the latter equality, which
is possible owing to integrability conditions, we deduce inversion
formula (3.13), where
$$\hat{{\cal M}}_{i\tau}(x)= {1\over \pi^3 i} \int_{b- i\infty}^{b+
i\infty} 2^{s-2}\Gamma\left({s+ i\tau\over 2}\right)\Gamma\left({s-
i\tau\over 2}\right){x^{-s}\over \zeta(-s)} ds, x >0.$$
To complete the proof, we will show that the kernel $\hat{{\cal
M}}_{i\tau}(x)$ can be written in the form (3.14). To do this we
appeal to the functional equation (1.12) for the Riemann
zeta-function and duplication formula for the Euler gamma -
function. Thus it gives
$$\hat{{\cal M}}_{i\tau}(x)=  - {1\over 2\pi^{5/2} i} \int_{b/2- i\infty}^{b/2+
i\infty} \frac{\Gamma\left(s+ {i\tau\over 2}\right)\Gamma\left(s-
{i\tau\over 2}\right)\Gamma(s) \Gamma\left(1-
s\right)}{\zeta(1+2s)\Gamma(s+ 1/2)\Gamma(1+s)}(x/2\pi)^{-2s}
ds.\eqno(3.21)$$
In the mean time the Parseval identity for the Mellin transform [11]
and relations (8.4.19.1), (8.4.23.27) in [5], Vol. 3 lead to the
equality
$$- {1\over 2\pi^{5/2} i} \int_{b/2- i\infty}^{b/2+
i\infty} \frac{\Gamma\left(s+ {i\tau\over 2}\right)\Gamma\left(s-
{i\tau\over 2}\right)\Gamma(s) \Gamma\left(1-
s\right)}{\zeta(1+2s)\Gamma(s+ 1/2)\Gamma(1+s)}(x/2\pi)^{-2s} ds$$
$$= -{2\over \pi^2}\int_0^\infty K^2_{i\tau/2}\left({xy\over
4\pi}\right)J_1(y)dy,$$
where $J_1(y)$ is the Bessel function of the first kind [1], Vol.
II. But the latter integral is calculated in [5], Vol. 2, relation
(2.16.43.2) and we obtain the result
$$ -{2\over \pi^2}\int_0^\infty K^2_{i\tau/2}\left({xy\over
4\pi}\right)J_1(y)dy= {2^{3/2}\over \sqrt x\
\tau\sinh(\pi\tau/2)}\left(1+ {4\pi^2\over x^2}\right)^{1/4}
P_{(i\tau-1)/2}^{1/2}\left(1+ {8\pi^2\over x^2}\right),\eqno(3.22)$$
where $P^d_\nu(z)$ is the associated Legendre function [1], Vol. I.
Hence, returning to (3.21) and combining with series (1.5), we
substitute it inside the integral. Then changing the order of
integration and summation via the absolute convergence and appealing
to (3.22), we come out with (3.14).
\end{proof}

\section{ Salem's type equivalences to the Riemann hypothesis}

In 1953 Salem [7] proved that the Riemann hypothesis is true, i.e.
the Riemann zeta-function $\zeta(s)$ is free of zeros in the strip
$1/2 < {\rm Re} s < 1$ is equivalent to the fact, that the
homogeneous integral equation
$$\int_0^\infty {y^{\delta-1}\over e^{xy} +1} h(y) dy=0, \ x >0, \ {1\over
2} < \delta < 1,\eqno(4.1)$$
has no nontrivial solutions in the space of bounded measurable
functions on $\mathbb{R_+}$. But after a simple change of variable
this equation becomes (2.15), where $g(x)=0$ and $f(t)=
t^{-\delta}h(1/t)$. Therefore reciprocities (2.9), (2.10) and
(2.15), (2.16) lead to

{\bf Corollary 1}. {\it Let $h(x)$ be a solution of homogeneous
equation $(4.1)$  such that $x^{-\delta}h(1/x) \in {\cal
M}^{-1}(L_c), \ c_0 > 1, \  1/2 < \delta < 1$. Then $h(x)\equiv 0$}.

\begin{proof} Indeed, there exists a function $h_\delta^*(s) \in
L_1(c )$ such that
$$x^{-\delta}h\left({1\over x}\right)= {1\over 2\pi i}\int_c
h_\delta^*(s) x^{-s} ds, \ x >0.$$
Hence
$$|h(x)| \le {x^{c_0-\delta}\over 2\pi} \int_c |h_\delta^*(s) ds|$$
and since $c_0 > \delta$, we have that $h(x)$ is continuous on
$\mathbb{R}_+$  and $h(x)=o(1), \ x \to 0$.  Applying inversion
formula (2.16) with $g=0$ we get the result.

\end{proof}

Let us prove the following equivalence to the Riemann hypothesis of
the Salem type.

{\bf Theorem 7}. {\it The Riemann hypothesis is true, if and only if
for any bounded measurable function $f(x)$ on $\mathbb{R}$
satisfying integral equation
$$ \int_{\mathbb{R}^2} \frac{e^{- \delta u}f(u)}{
(e^{e^{x-u-t}}+1)(e^{e^t}+1)}du dt =0, {1\over 2} < \delta <
1,\eqno(4.2)$$
for all $x \in \mathbb{R}$ it follows that $f$ is zero almost
everywhere.}

\begin{proof} Calling again (1.14) and properties of the Mellin
transform and its convolution [11], [14] it is not difficult to
derive the equality
$$[(1-2^{1-s})\zeta(s)\Gamma(s)]^2= \int_0^\infty t^{s-1}\\
\int_0^\infty {du\over  u (e^{t/u}+1)(e^{u}+1)} dt, \ {\rm Re} \ s
>0.\eqno(4.3)$$
On the other hand, the reciprocal inversion of the Mellin transform
yields
$$ \int_0^\infty {du\over  u (e^{x/u}+1)(e^{u}+1)}\\
= {1\over 2\pi i} \int_{\mu- i\infty}^{\mu+
i\infty}[(1-2^{1-s})\zeta(s)\Gamma(s)]^2 x^{-s}ds.\eqno(4.4)$$
The left-hand side of (4.4) is positive and via (4.3)
\begin{equation*}
\begin{split}
\int_0^\infty \int_0^\infty {t^{\delta-1} du dt\over  u
(e^{t/u}+1)(e^{u}+1)} =
[(1-2^{1-\delta})\zeta(\delta)\Gamma(\delta)]^2,
\end{split}
\end{equation*}
which after a simple change of variables is equivalent to the
condition
$$\int_{\mathbb{R}^2} {e^{\delta y} du dy \over (e^{e^{y-u}}+1)(e^{e^u}+1)} < \infty.$$
Hence following  as in [7]  Wiener's ideas  about an equivalence of
the completeness in $L_1(\mathbb{R})$ of translations
$$ e^{\delta (x-y)} \int_{\mathbb{R}} {du\over (e^{e^{x-y-u}}+1)(e^{e^u}+1)}, \ x \in \mathbb{R} $$
and the absence of zeros of $[(1-2^{1-s})\zeta(s)\Gamma(s)]^2$, i.e.
zeros of $\zeta(s)$ in the critical strip $1/2 < {\rm Re} \ s < 1$,
we complete the proof.

\end{proof}

{\bf Remark 1}. Reminding integral representation (3.3) of the
modified Bessel function and    invoking identity (1.4),  one can we
write equality (4.4) in the form
$$
{1\over 2} \int_0^\infty {du\over  u (e^{x/u}+1)(e^{u}+1)}=
\sum_{n=1}^\infty d(n) \left[K_0(2\sqrt {nx}) \right. \\ \left. - 4
K_0(2\sqrt {2nx})+ 4 K_0(4\sqrt {nx})\right].\eqno(4.5)$$

Hence substituting (4.5) into (4.2), we change the order of
integration and summation via absolute and uniform convergence since
(see Section 1) $d(n)= O(n^\varepsilon), \varepsilon
> 0, n \to \infty$.  Consequently,  Theorem 7 can be reformulated as

{\bf Theorem 8}. {\it  The Riemann hypothesis is true, if and only
if for any bounded measurable function $f(x)$ on $\mathbb{R}$ and
all $x \in \mathbb{R}$ the equation
\begin{equation*}
\sum_{n=1}^\infty d(n) \left[({\cal K}_n f)(x)- 4 ({\cal K}_{2n}
f)(x)+ 4 ({\cal K}_{4n} f)(x)\right]=0,
\end{equation*}
where
\begin{equation*}
\begin{split}
({\cal K}_n f)(x)= \int_{-\infty}^\infty e^{- \delta u}
K_0\left(2\sqrt {n}\  e^{(x-u)/2}\right)f(u)du,\  {1\over 2} <
\delta < 1,
\end{split}
\end{equation*}
is the Meijer type convolution transform $[3]$, has no nontrivial
solutions.}

Finally a class of Salem's type equivalences to the Riemann
hypothesis is given by

{\bf Theorem 9}. {\it  Let $k, m \in \mathbb{N}_0, \ k \le m$ and
the kernel $U_{k,m}(x), \ x >0$ is defined by formulas $(2.19),
(2.20)$, correspondingly. The Riemann hypothesis is true, if and
only if for any bounded measurable function $f(x)$ on $\mathbb{R}$
satisfying integral equation
$$\int_{\mathbb{R}} e^{- \delta u}U_{k,m}\left(e^{x-u}\right)f(u)du
=0,\  {1\over 2} < \delta < 1,\eqno(4.6)$$
 for all $x \in \mathbb{R}$ it follows that $f$ is zero almost
 everywhere.}

\begin{proof} Employing  inversion formula (1.15) of  the Mellin
transform, we derive, reciprocally, from (2.17)
\begin{equation*}
\begin{split}
[(1-2^{1-s})\zeta(s)]^{k+1}\Gamma^{m+1}(s) = \int_0^\infty
U_{k,m}(t) t^{s-1} dt,\  {\rm Re} \ s >0.
\end{split}
\end{equation*}
Moreover, $U_{k,m}(x),\ x >0$ is positive (see  (2.19), (2.20)) and
for $\delta \in (1/2,1)$
\begin{equation*}
\begin{split}
\int_0^\infty U_{k,m}(t) t^{\sigma-1} dt  =
[(1-2^{1-\sigma})\zeta(\sigma)]^{k+1}\Gamma^{m+1}(\sigma).
\end{split}
\end{equation*}
This yields
$$\int_{\mathbb{R}} e^{\delta y} U_{k,m}(e^y) dy  < \infty.$$
Hence as in Theorem 6 the completeness in $L_1(\mathbb{R})$ of
translations
$$ e^{\delta (x-y)} U_{k,m}(e^{x-y}), \ x \in \mathbb{R} $$
is equivalent to the absence of zeros of
$[(1-2^{1-s})\zeta(s)]^{k+1}\Gamma^{m+1}(s)$, i.e. zeros of
$\zeta(s)$ in the critical strip $1/2 < \delta  < 1$.

\end{proof}

\bigskip
\centerline{{\bf Acknowledgments}}
\bigskip
The present investigation was supported, in part,  by the "Centro de
Matem{\'a}tica" of the University of Porto.

\bigskip
\centerline{{\bf References}}
\bigskip
\baselineskip=12pt
\medskip
\begin{enumerate}

\item[{\bf 1.}\ ]
 A. Erd\'elyi, W. Magnus, F. Oberhettinger and F.G. Tricomi,
{\it Higher Transcendental Functions}, Vols. I and  II, McGraw-Hill,
New York, London and Toronto (1953).

\item[{\bf 2.}\ ] R.R. Goldberg, Inversions of generalized Lambert transforms,
 {\it Duke Math. J.}, {\bf 25} (1958), 459- 476.

\item[{\bf 3.}\ ]I.I. Hirschman and D.V. Widder, {\it The Convolution Transform}, Princeton University Press, Princeton, New Jersey (1955).

\item[{\bf 4.}\ ] W.B. Pennington, Widder's inversion formula for
the Lambert transform,  {\it Duke Math. J.}, {\bf 27} (1960),
561-568.

\item[{\bf 5.}\ ]  A.P. Prudnikov, Yu. A. Brychkov and O. I.
Marichev, {\it Integrals and Series: Vol. 2: Special Functions},
Gordon and Breach, New York  (1986); {\it Vol. 3: More Special
Functions}, Gordon and Breach, New York  (1990).

\item[{\bf 6.}\ ] S. Ramanujan, Some formulas in the analytic theory of numbers,
{\it Messenger of  Math.}, {\bf 45} (1916), 81-84.

\item[{\bf 7.}\ ] R. Salem, Sur une proposition equivalente a l'hipothese de Riemann.
{\it Comptes Rendus Math.}, {\bf 236} (1953), 1127-1128.

\item[{\bf 8.}\ ] J. Sandor, D.S. Mitrinovich, B. Crstici, {\it
Handbook of Number Theory}, Vol. I, Springer, Dordrecht (2006).

\item[{\bf 9.}\ ]I.N. Sneddon, {\it The Use of Integral
Transforms}, McGray Hill, New York (1972).

\item[{\bf 10.}\ ] E.C. Titchmarsh, {\it The Theory of The Riemann Zeta- Function},
The Clarendon Press, Oxford, Second edition  (1986).

\item[{\bf 11.}\ ]  E.C. Titchmarsh, {\it  An Introduction to the
Theory of Fourier Integrals}, Clarendon Press, Oxford ( 1937).

\item[{\bf 12.}\ ] Vu Kim Tuan, O.I. Marichev and S.B. Yakubovich, Composition structure of integral
transformations,  {\it J. Soviet Math.}, {\bf 33} (1986), 166-169.

\item[{\bf 13.}\ ] D. V. Widder, An inversion of the Lambert transform,
 {\it Math. Mag.}, {\bf 23} (1950), 171-182.

\item[{\bf 14.}\ ] S. B. Yakubovich and  Yu. F. Luchko, {\it The
Hypergeometric Approach to Integral Transforms and Convolutions}.
Mathematics and its Applications, 287. Kluwer Academic Publishers
Group, Dordrecht (1994).

\item[{\bf 15.}\ ] S.B. Yakubovich and B. Fisher, A class of index transforms with general kernels,
 {\it Math. Nachr.}, {\bf 200} (1999), 165-182.

\item[{\bf 16.}\ ] S.B. Yakubovich, {\it Index Transforms},  World
Scientific Publishing Company, Singapore, New Jersey, London and
Hong Kong (1996).

\end{enumerate}

\vspace{5mm}

\noindent S.Yakubovich\\
Department of Mathematics,\\
Faculty of Sciences,\\
University of Porto,\\
Campo Alegre st., 687\\
4169-007 Porto\\
Portugal\\
E-Mail: syakubov@fc.up.pt\\

\end{document}